# LASSO-TYPE RECOVERY OF SPARSE REPRESENTATIONS FOR HIGH-DIMENSIONAL DATA

BY NICOLAI MEINSHAUSEN[1] AND BIN YU[2]

*University of Oxford and University of California, Berkeley*

The Lasso is an attractive technique for regularization and variable selection for high-dimensional data, where the number of predictor variables $p_n$ is potentially much larger than the number of samples $n$. However, it was recently discovered that the sparsity pattern of the Lasso estimator can only be asymptotically identical to the true sparsity pattern if the design matrix satisfies the so-called *irrepresentable condition*. The latter condition can easily be violated in the presence of highly correlated variables.

Here we examine the behavior of the Lasso estimators if the *irrepresentable condition* is relaxed. Even though the Lasso cannot recover the correct sparsity pattern, we show that the estimator is still consistent in the $\ell_2$-norm sense for fixed designs under conditions on (a) the number $s_n$ of nonzero components of the vector $\beta_n$ and (b) the minimal singular values of design matrices that are induced by selecting small subsets of variables. Furthermore, a rate of convergence result is obtained on the $\ell_2$ error with an appropriate choice of the smoothing parameter. The rate is shown to be optimal under the condition of bounded maximal and minimal sparse eigenvalues. Our results imply that, with high probability, all important variables are selected. The set of selected variables is a meaningful reduction on the original set of variables. Finally, our results are illustrated with the detection of closely adjacent frequencies, a problem encountered in astrophysics.

**1. Introduction.** The Lasso was introduced by [29] and has since been proven to be very popular and well studied [18, 35, 41, 42]. Some reasons for the popularity might be that the entire regularization path of the Lasso can be computed efficiently [11, 25], that Lasso is able to handle more predictor

---

Received December 2006; revised December 2007.
[1]Supported by DFG (Deutsche Forschungsgemeinschaft).
[2]Supported in part by a Guggenheim fellowship and Grants NSF DMS-06-05165 (06-08), NSF DMS-03-036508 (03-05) and ARO W911NF-05-1-0104 (05-07).

variables than samples and produces sparse models which are easy to interpret. Several extensions and variations have been proposed [5, 21, 36, 40, 42].

1.1. *Lasso-type estimation.* The Lasso estimator, as introduced by [29], is given by

$$\hat{\beta}^\lambda = \arg\min_{\beta} \|Y - X\beta\|_{\ell_2}^2 + \lambda\|\beta\|_{\ell_1}, \tag{1}$$

where $X = (X_1, \ldots, X_p)$ is the $n \times p$ matrix whose columns consist of the $n$-dimensional fixed predictor variables $X_k$, $k = 1, \ldots, p$. The vector $Y$ contains the $n$-dimensional set of real-valued observations of the response variable.

The distribution of Lasso-type estimators has been studied in Knight and Fu [18]. Variable selection and prediction properties of the Lasso have been studied extensively for high-dimensional data with $p_n \gg n$, a frequently encountered challenge in modern statistical applications. Some studies Bunea, Tsybakov and Wegkamp, for example, [2], Greenshtein and Ritov, for example, [13], van de Geer, for example, [34] have focused mainly on the behavior of prediction loss. Much recent work aims at understanding the Lasso estimates from the point of view of model selection, including Candes and Tao [5], Donoho, Elad and Temlyakov [10], Meinshausen and Bühlmann [23], Tropp [30], Wainwright [35], Zhao and Yu [41], Zou [42]. For the Lasso estimates to be close to the model selection estimates when the data dimensions grow, all the aforementioned papers assumed a sparse model and used various conditions that require the irrelevant variables to be not too correlated with the relevant ones. Incoherence is the terminology used in the deterministic setting of Donoho, Elad and Temlyakov [10] and "irrepresentability" is used in the stochastic setting (linear model) of Zhao and Yu [41]. Here we focus exclusively on the properties of the estimate of the coefficient vector under squared error loss and try to understand the behavior of the estimate under a relaxed *irrepresentable condition* (hence we are in the stochastic or linear model setting). The aim is to see whether the Lasso still gives meaningful models in this case.

More discussions on the connections with other works will be covered in Section 1.5 after notions are introduced to state explicitly what the irrepresentable condition is so that the discussions are clearer.

1.2. *Linear regression model.* We assume a linear model for the observations of the response variable $Y = (Y_1, \ldots, Y_n)^T$,

$$Y = X\beta + \varepsilon, \tag{2}$$

where $\varepsilon = (\varepsilon_1, \ldots, \varepsilon_n)^T$ is a vector containing independently and identically distributed noise with $\varepsilon_i \sim \mathcal{N}(0, \sigma^2)$ for all $i = 1, \ldots, n$. The assumption of Gaussianity could be relaxed and replaced with exponential tail bounds on the noise if, additionally, predictor variables are assumed to be bounded. When there is a question of nonidentifiability for $\beta$, for $p_n > n$, we define $\beta$ as

$$\beta = \underset{\{\beta\,:\,EY=X\beta\}}{\arg\min} \|\beta\|_{\ell_1}. \tag{3}$$

The aim is to recover the vector $\beta$ as well as possible from noisy observations $Y$. For the equivalence between $\ell_1$- and $\ell_0$-sparse solutions see, for example, Donoho [8], Donoho and Elad [9], Fuchs [12], Gribonval and Nielsen [14], Tropp [30, 31].

1.3. *Recovery of the sparsity pattern and the irrepresentable condition.* There is empirical evidence that many signals in high-dimensional spaces allow for a sparse representation. As an example, wavelet coefficients of images often exhibit exponential decay, and a relatively small subset of all wavelet coefficients allow a good approximation to the original image [17, 19, 20]. For conceptual simplicity, we assume in our regression setting that the vector $\beta$ is sparse in the $\ell_0$-sense and many coefficients of $\beta$ are identically zero. The corresponding variables have thus no influence on the response variable and could be safely removed. The sparsity pattern of $\beta$ is understood to be the sign function of its entries, with $\mathrm{sign}(x) = 0$ if $x = 0$, $\mathrm{sign}(x) = 1$ if $x > 0$ and $\mathrm{sign}(x) = -1$ if $x < 0$. The sparsity pattern of a vector might thus look like

$$\mathrm{sign}(\beta) = (+1, -1,\ 0,\ 0, +1, +1, -1, +1,\ 0,\ 0, \ldots),$$

distinguishing whether variables have a positive, negative or no influence at all on the response variable. It is of interest whether the sparsity pattern of the Lasso estimator is a good approximation to the true sparsity pattern. If these sparsity patterns agree asymptotically, the estimator is said to be *sign consistent* [41].

DEFINITION 1 (Sign consistency). An estimator $\hat{\beta}^\lambda$ is *sign consistent* if and only if

$$P\{\mathrm{sign}(\beta) = \mathrm{sign}(\hat{\beta})\} \to 1 \qquad \text{as } n \to \infty.$$

It was shown independently in Zhao and Yu [41] and Zou [42] in the linear model case and [23] in a Gaussian Graphical Model setting that *sign consistency* requires a condition on the design matrix. The assumption was termed *neighborhood stability* in Meinshausen and Bühlmann [23] and *irrepresentable condition* in Zhao and Yu [41]. Let $C = n^{-1} X^T X$. The dependence on $n$ is neglected notationally.

4         N. MEINSHAUSEN AND B. YU

DEFINITION 2 (Irrepresentable condition). Let $K = \{k : \beta_k \neq 0\}$ be the set of relevant variables and let $N = \{1, \ldots, p\} \setminus K$ be the set of noise variables. The sub-matrix $C_{HK}$ is understood as the matrix obtained from $C$ by keeping rows with index in the set $H$ and columns with index in $K$. The *irrepresentable condition* is fulfilled if

$$\|C_{NK} C_{KK}^{-1} \operatorname{sign}(\beta_K)\|_{\ell_\infty} < 1.$$

In Zhao and Yu [41], an additional *strong irrepresentable condition* is defined which requires that the above elements are not merely smaller than 1 but are uniformly bounded away from 1. Zhao and Yu [41], Zou [42] and Meinshausen and Bühlmann [23] show that the Lasso is sign consistent only if the *irrepresentable condition* holds.

PROPOSITION 1 (Sign consistency).  *Assume that the irrepresentable condition or neighborhood stability is not fulfilled. Then there exists no sequence $\lambda = \lambda_n$ such that the estimator $\hat\beta^\lambda$ is sign consistent.*

It is worth noting that a slightly stronger condition has been used in Tropp [30, 31] in a deterministic study of Lasso's model selection properties where $1 - C_{NK} C_{KK}^{-1}$ is called ERC (exact recovery coefficient). A positive ERC implies the irrepresentable condition for all $\beta$ values.

In practice, it might be difficult to verify whether the condition is fulfilled. This led various authors to propose interesting extensions to the Lasso [22, 39, 42]. Before giving up on the Lasso altogether, however, we want to examine in this paper in what sense the original Lasso procedure still gives sensible results, even if the *irrepresentable condition* or, equivalently, *neighborhood stability* is not fulfilled.

1.4. *$\ell_2$-consistency.* The aforementioned studies showed that if the *irrepresentable condition* is not fulfilled, the Lasso cannot select the correct sparsity pattern. In this paper we show that the Lasso selects in these cases the nonzero entries of $\beta$ *and some not-too-many additional* zero entries of $\beta$ under relaxed conditions than the irrepresentable condition. The nonzero entries of $\beta$ are in any case included in the selected model. Moreover, the size of the estimated coefficients allows to separate the few truly zero and the many nonzero coefficients. However, we note that in extreme cases, when the variables are linearly dependent, even these relaxed conditions will be violated. In these situations, it is not sensible to use the $\ell_2$-metric on $\beta$ to assess Lasso.

Our main result shows the $\ell_2$-consistency of the Lasso, even if the *irrepresentable condition* is violated. To be precise, an estimator is said to be $\ell_2$-consistent if

(4)           $\|\hat\beta - \beta\|_{\ell_2} \to 0 \qquad \text{as } n \to \infty.$

Rates of convergence results will also be derived and under the condition of bounded maximal and minimal sparse eigenvalues, the rate is seen optimal. An $\ell_2$-consistent estimator is attractive, as important variables are chosen with high probability and falsely chosen variables have very small coefficients. The bottom line will be that even if the sparsity pattern of $\beta$ cannot be recovered by the Lasso, we can still obtain a good approximation.

1.5. *Related work.* Prediction loss for high-dimensional regression under an $\ell_1$-penalty has been studied for quadratic loss function in Greenshtein and Ritov [13] and for general Lipschitz loss functions in van de Geer [34]. With a focus on aggregation, similarly interesting results are derived in Bunea, Tsybakov and Wegkamp [3]. Both van de Geer [34] and Bunea, Tsybakov and Wegkamp [3] obtain impressive results for random design and sharp bounds for the $\ell_1$-distance between the vector $\beta$ and its Lasso estimate $\hat{\beta}^\lambda$. In the current manuscript, we focus on the $\ell_2$-estimation loss on $\beta$. As a consequence, we can derive consistency in the sense of (4) under the condition that $s_n \log p_n/n \to 0$ for $n \to \infty$ (ignoring $\log n$ factors). An implication of our work is thus that the sparsity $s_n$ is allowed to grow almost as fast as the sample size if one is interested to obtain convergence in $\ell_2$-norm. In contrast, the results in [3, 34] require $s_n = o(\sqrt{n})$ to obtain convergence in $\ell_1$-norm.

The recent independent work of Zhang and Huang [38] shows that the subspace spanned by the variables selected by Lasso is close to an optimal subspace. The results also imply that important variables are chosen with high probability and provides a tight bound on the $\ell_2$-distance between the vector $\beta$ and its Lasso estimator. A "partial Riesz condition" is employed in [38], which is rather similar to our notion of *incoherent design*, defined further below in (6).

We would like to compare the results of this manuscript briefly with results in Donoho [8] and Candes and Tao [5], as both of these papers derive bounds on the $\ell_2$-norm distance between $\beta$ and $\hat{\beta}$ for $\ell_1$-norm constrained estimators. In Donoho [8] the design is random and the random predictor variables are assumed to be independent. The results are thus not directly comparable to the results derived here for general fixed designs. Nevertheless, results in Meinshausen and Bühlmann [23] suggest that the *irrepresentable condition* is with high probability fulfilled for independently normal distributed predictor variables. The results in Donoho [8] can thus not directly be used to study the behavior of the Lasso under a violated *irrepresentable condition*, which is our goal in the current manuscript.

Candes and Tao [5] study the properties of the so-called "Dantzig selector," which is very similar to the Lasso, and derive bounds on the $\ell_2$-distance between the vector $\beta$ and the proposed estimator $\hat{\beta}$. The results are derived under the condition of a *Uniform Uncertainty Principle* (UUP), which was introduced in Candes and Tao [4]. The UUP is related to our assumptions on sparse eigenvalues in this manuscript. A comparison between these two assumptions is given after the formulation (10) of the UUP. The bounds on the $\ell_2$-distance between the true coefficient vector $\beta$ and its Lasso estimator (obtained in the current manuscript) or, respectively, "Dantzig selector" (obtained in [5]) are quite similar in nature. This comes maybe as no surprise since the formulation of the "Dantzig selector" is quite similar to the Lasso [24]. However, it does not seem straightforward to translate the bounds obtained for the "Dantzig selector" into bounds for the Lasso estimator and vice versa. We employ also somewhat different conditions because there could be situations of design matrix arising in statistical practice where the dependence between the predictors is stronger than what is allowed by the UUP, but would satisfy our condition of "incoherent design" to be defined in the next section. It would certainly be of interest to study the connection between the Lasso and "Dantzig selector" further, as the solutions share many similarities.

Final note: a recent follow-up work [1] provides similar bounds as in this paper for both Lasso and Dantzig selectors.

**2. Main assumptions and results.** First, we introduce the notion of *sparse eigenvalues*, which will play a crucial role in providing bounds for the convergence rates of the Lasso estimator. Thereafter, the assumptions are explained in detail and the main results are given.

2.1. *Sparse eigenvalues.* The notion of *sparse eigenvalues* is not new and has been used before [8]; we merely intend to fixate notation. The *m-sparse minimal eigenvalue* of a matrix is the minimal eigenvalue of any $m \times m$-dimensional submatrix.

DEFINITION 3. The *m-sparse minimal eigenvalue* and *m-sparse maximal eigenvalue* of $C$ are defined as

$$(5) \quad \phi_{\min}(m) = \min_{\beta \,:\, \|\beta\|_{\ell_0} \leq \lceil m \rceil} \frac{\beta^T C \beta}{\beta^T \beta} \quad \text{and} \quad \phi_{\max}(m) = \max_{\beta \,:\, \|\beta\|_{\ell_0} \leq \lceil m \rceil} \frac{\beta^T C \beta}{\beta^T \beta}.$$

The minimal eigenvalue of the unrestricted matrix $C$ is equivalent to $\phi_{\min}(p)$. If the number of predictor variables $p_n$ is larger than sample size, $p_n > n$, this eigenvalue is zero, as $\phi_{\min}(m) = 0$ for any $m > n$.

A crucial factor contributing to the convergence of the Lasso estimator is the behavior of the smallest *m-sparse eigenvalue*, where the number $m$ of variables over which the minimal eigenvalues is computed is roughly the same order as the sparsity $s_n$, or the number of nonzero components, of the true underlying vector $\beta$.

2.2. *Sparsity multipliers and incoherent designs.* As apparent from the interesting discussion in Candes and Tao [5], one cannot allow arbitrarily large "coherence" between variables if one still hopes to recover the correct sparsity pattern. Assume that there are two vectors $\beta$ and $\tilde{\beta}$ so that the signal can be represented by either vector $X\beta = X\tilde{\beta}$ and both vectors are equally sparse, say $\|\beta\|_{\ell_0} = \|\tilde{\beta}\|_0 = s_n$ and are not identical. We have no hope of distinguishing between $\beta$ and $\tilde{\beta}$ in such a case: if indeed $X\beta = X\tilde{\beta}$ and $\beta$ and $\tilde{\beta}$ are not identical, it follows that the minimal sparse eigenvalue $\phi_{\min}(2s_n) = 0$ vanishes as $X(\beta - \tilde{\beta}) = 0$ and $\|\beta - \tilde{\beta}\|_{\ell_0} \leq 2s_n$. If the minimal sparse eigenvalue of a selection of $2s_n$ variables is zero, we have no hope of recovering the true sparse underlying vector from noisy observations.

To define our assumption about sufficient conditions for recovery, we need the definition of *incoherent design*. As motivated by the example above, we would need a lower bound on the minimal eigenvalue of at least $2s_n$ variables, where $s_n$ is again the number of nonzero coefficients. We now introduce the concepts of *sparsity multiplier* ad *incoherent design* to make this requirement a bit more general, as minimal eigenvalues are allowed to converge to zero slowly.

A design is called *incoherent* in the following if minimal sparse eigenvalues are not decaying too fast, in a sense made precise in the definition below. For notational simplicity, let in the following

$$\phi_{\max} = \phi_{\max}(s_n + \min\{n, p_n\})$$

be the maximal eigenvalue of a selection of at most $s_n + \min\{n, p_n\}$ variables. At the cost of more involved proofs, one could also work with the maximal eigenvalue of a smaller selection of variables instead. Even though we do not assume an upper bound for the quantity $\phi_{\max}$, it would not be very restrictive to do so for the $p_n \gg n$ setting. To be specific, assume multivariate normal predictors. If the maximal eigenvalue of the population covariance matrix, which is induced by selecting $2n$ variables, is bounded from above by an arbitrarily large constant, it follows by Theorem 2.13 in Davidson and Szarek [7] or Lemma A3.1 in Paul [26] that the condition number of the induced sample covariance matrix observes a Gaussian tail bound. Using an entropy bound for the possible number of subsets when choosing $n$ out of $p_n$ variables. The maximal eigenvalue of a selection of $2\min\{n, p\}$ variables is thus bounded from above by some constant, with probability converging to 1 for $n \to \infty$ under the condition that $\log p_n = o(n^\kappa)$ for some $\kappa < 1$, and the assumption of a bounded $\phi_{\max}$, even though not needed, is thus maybe not overly restrictive.

As the maximal sparse eigenvalue is typically growing only very slowly as a function of the number of variables, the focus will be on the decay of the smallest sparse eigenvalue, which is a much more pressing problem for high-dimensional data.

DEFINITION 4 (Incoherent designs). A design is called *incoherent* if there exists a positive sequence $e_n$, the so-called *sparsity multiplier* sequence, such that

$$\text{(6)} \qquad \liminf_{n \to \infty} \frac{e_n \phi_{\min}(e_n^2 s_n)}{\phi_{\max}(s_n + \min\{n, p_n\})} \geq 18.$$

Our main result will require *incoherent design*. The constant 18 could quite possibly be improved upon. We will assume for the following that the multiplier sequence is the smallest. Below, we give some simple examples under which the condition of *incoherent design* is fulfilled.

2.2.1. *Example: block designs.* The first example is maybe not overly realistic but gives, hopefully, some intuition for the condition. A "block design" is understood to have the structure

$$\text{(7)} \qquad n^{-1} X^T X = \begin{pmatrix} \Sigma(1) & 0 & \cdots & 0 \\ 0 & \Sigma(2) & \cdots & 0 \\ \cdots & \cdots & \cdots & \cdots \\ 0 & 0 & \cdots & \Sigma(d) \end{pmatrix},$$

where the matrices $\Sigma(1), \ldots, \Sigma(d)$ are of dimension $b(1), \ldots, b(d)$, respectively. The minimal and maximal eigenvalues over all $d$ sub-matrices are denoted by

$$\phi_{\min}^{\text{block}} := \min_k \min_{u \in \mathbb{R}^{b(k)}} \frac{u^T \Sigma(k) u}{u^T u}, \qquad \phi_{\max}^{\text{block}} := \max_k \max_{u \in \mathbb{R}^{b(k)}} \frac{u^T \Sigma(k) u}{u^T u}.$$

In our setup, all constants are allowed to depend on the sample size $n$. The question arises if simple bounds can be found under which the design is *incoherent* in the sense of (6). The blocked sparse eigenvalues are trivial lower and upper bounds, respectively, for $\phi_{\min}(u)$ and $\phi_{\max}(u)$ for all values of $u$. Choosing $e_n$ such that $e_n^2 s_n = o(n)$, the condition (6) of *incoherent design* requires then $e_n \phi_{\min}(e_n^2 s_n) \gg \phi_{\max}(s_n + \min\{n, p_n\})$. Using $\phi_{\min}(e_n^2 s_n) \geq \phi_{\min}^{\text{block}}$ and $\phi_{\max} \leq \phi_{\max}^{\text{block}}$, it is sufficient if there exists a sequence $e_n$ with $e_n = o(\phi_{\max}^{\text{block}}/\phi_{\min}^{\text{block}})$. Together with the requirement $e_n^2 s_n = o(n)$, the condition of *incoherent design* is fulfilled if, for $n \to \infty$,

$$\text{(8)} \qquad s_n = o\left(\frac{n}{c_n^2}\right),$$

where the condition number $c_n$ is given by

$$\text{(9)} \qquad c_n := \phi_{\max}^{\text{block}}/\phi_{\min}^{\text{block}}.$$

Under increasingly stronger assumption on the sparsity, the condition number $c_n$ can thus grow almost as fast as $\sqrt{n}$, while still allowing for *incoherent design*.

2.2.2. *More examples of incoherent designs.* Consider two more examples of incoherent design:

- The condition (6) of *incoherent design* is fulfilled if the minimal eigenvalue of a selection of $s_n(\log n)^2$ variables is vanishing slowly for $n \to \infty$ so that

$$\phi_{\min}\{s_n(\log n)^2\} \gg \frac{1}{\log n} \phi_{\max}(s_n + \min\{p_n, n\}).$$

- The condition is also fulfilled if the minimal eigenvalue of a selection of $n^\alpha s_n$ variables is vanishing slowly for $n \to \infty$ so that

$$\phi_{\min}(n^\alpha s_n) \gg n^{-\alpha/2} \phi_{\max}.$$

These results can be derived from (6) by choosing the sparse multiplier sequences $e_n = \log n$ and $e_n = n^{\alpha/2}$, respectively. Some more scenarios of *incoherent design* can be seen to satisfy (6).

2.2.3. *Comparison with the uniform uncertainty principle.* Candes and Tao [5] use a *Uniform Uncertainty Principle* (UUP) to discuss the convergence of the so-called *Dantzig selector*. The UUP can only be fulfilled if the minimal eigenvalue of a selection of $s_n$ variables is bounded from below by a constant, where $s_n$ is again the number of nonzero coefficients of $\beta$. In the original version, a necessary condition for UUP is

(10) $$\phi_{\min}(s_n) + \phi_{\min}(2s_n) + \phi_{\min}(3s_n) > 2.$$

At the same time, a bound on the maximal eigenvalue is a condition for the UUP in [5],

(11) $$\phi_{\max}(s_n) + \phi_{\max}(2s_n) + \phi_{\max}(3s_n) < 4.$$

This UUP condition is different from our incoherent design condition. In some sense, the UUP is weaker than *incoherent design*, as the minimal eigenvalues are calculated over only $3s_n$ variables. In another sense, UUP is quite strong as it demands, in form (10) and assuming $s_n \geq 2$, that *all pairwise* correlations between variables be less than 1/3! The condition of *incoherent design* is weaker as the eigenvalue can be bounded from below by an arbitrarily small constant (as opposed to the large value implied by the UUP). Sparse eigenvalues can even converge slowly to zero in our setting.

Taking the example of block designs from further above, *incoherent design* allowed for the condition number (9) to grow almost as fast as $\sqrt{n}$. In contrast, if the sparsity $s_n$ is larger than the maximal block-size, the UUP requires that the condition number $c_n$ be bounded from above by a positive constant. Using its form (10) and the corresponding bound (11) for the maximal eigenvalue, it implies specifically that $c_n \leq 2$, which is clearly stricter than the condition (8).

2.2.4. *Incoherent designs and the irrepresentable condition.* One might ask in what sense the notion of *incoherent design* is more general than the *irrepresentable condition*. At first, it might seem like we are simply replacing the strict condition of *irrepresentable condition* by a similarly strong condition on the design matrix.

Consider first the classical case of a fixed number $p_n$ of variables. If the covariance matrix $C = C_n$ is converging to a positive definite matrix for large sample sizes, the design is automatically *incoherent*. On the other hand, it is easy to violate the *irrepresentable condition* in this case; for examples, see Zou [42].

The notion of *incoherent* designs is only a real restriction in the high-dimensional case with $p_n > n$. Even then, it is clear that the notion of *incoherence* is a relaxation from *irrepresentable condition*, as the *irrepresentable condition* can easily be violated even though all sparse eigenvalues are bounded well away from zero.

2.3. *Main result for high-dimensional data ($p_n > n$).* We first state our main result.

THEOREM 1 (Convergence in $\ell_2$-norm). *Assume the incoherent design condition (6) with a sparsity multiplier sequence $e_n$. If $\lambda \propto \sigma e_n \sqrt{n \log p_n}$, there exists a constant $M > 0$ such that, with probability converging to 1 for $n \to \infty$,*

$$\|\beta - \hat{\beta}^{\lambda_n}\|_{\ell_2}^2 \leq M \sigma^2 \frac{s_n \log p_n}{n} \frac{e_n^2}{\phi_{\min}^2(e_n^2 s_n)}. \tag{12}$$

A proof is given in Section 3. It can be seen from the proofs that nonasymptotic bounds could be obtained with essentially the same results.

If we choose the smallest possible multiplier sequence $e_n$, one obtains not only the required lower bound $e_n \geq 18 \phi_{\max}/\phi_{\min}(e_n^2 s_n)$ from (6) but also an upper bound $e_n \leq K \phi_{\max}/\phi_{\min}(e_n^2 s_n)$. Plugging this into (12) yields the probabilistic bound, for some positive $M$,

$$\|\beta - \hat{\beta}^{\lambda_n}\|_{\ell_2}^2 \leq M \sigma^2 \frac{s_n \log p_n}{n} \frac{\phi_{\max}^2}{\phi_{\min}^4(e_n^2 s_n)}.$$

It is now easy to see that the convergence rate is essentially optimal as long as the relevant eigenvalues are bounded.

COROLLARY 1. *Assume that there exist constants $0 < \kappa_{\min} \leq \kappa_{\max} < \infty$ such that*

$$\liminf_{n \to \infty} \phi_{\min}(s_n \log n) \geq \kappa_{\min} \quad \text{and}$$
$$\limsup_{n \to \infty} \phi_{\max}(s_n + \min\{n, p_n\}) \leq \kappa_{\max}. \tag{13}$$

Then, for $\lambda \propto \sigma\sqrt{n\log p_n}$, there exists a constant $M > 0$ such that, with probability converging to 1 for $n \to \infty$,

$$\|\beta - \hat{\beta}^{\lambda_n}\|^2_{\ell_2} \leq M\sigma^2 \frac{s_n \log p_n}{n}.$$

The proof of this follows from Theorem 1 by choosing a constant *sparsity multiplier* sequence, for example, $20\kappa_{\max}/\kappa_{\min}$.

The rate of convergence achieved is essentially optimal. Ignoring the $\log p_n$ factor, it corresponds to the rate that could be achieved with maximum likelihood estimation if the true underlying sparse model would be known.

It is perhaps also worthwhile to make a remark about the penalty parameter sequence $\lambda$ and its, maybe unusual, reliance on the sparsity multiplier sequence $e_n$. If both the relevant minimal and maximal sparse eigenvalues in (6) are bounded from below and above, as in Corollary 1 above, the sequence $e_n$ is simply a constant. Any deviation from the usually optimal sequence $\lambda \propto \sigma\sqrt{n\log p_n}$ occurs thus only if the minimal sparse eigenvalues are decaying to zero for $n \to \infty$, in which case the penalty parameter is increased slightly. The value of $\lambda$ can be computed, in theory, without knowledge about the true $\beta$. Doing so in practice would not be a trivial task, however, as the sparse eigenvalues would have to be known. Moreover, the noise level $\sigma$ would have to be estimated from data, a difficult task for high-dimensional data with $p_n > n$. From a practical perspective, we mostly see the results as implying that the $\ell_2$-distance can be small for *some* value of the penalty parameter $\lambda$ along the solution path.

2.4. *Number of selected variables.* As a result of separate interest, it is perhaps noteworthy that bounds on the number of selected variables are derived for the proof of Theorem 1. For the setting of Corollary 1 above, where a constant *sparsity multiplier* can be chosen, Lemma 5 implies that, with high probability, at most $O(s_n)$ variables are selected by the Lasso estimator. The selected subset is hence of the same order of magnitude as the set of "truly nonzero" coefficients. In general, with high probability, no more than $e_n^2 s_n$ variables are selected.

2.5. *Sign consistency with two-step procedures.* It follows from our results above that the Lasso estimator can be modified to be *sign consistent* in a two-step procedure even if the *irrepresentable condition* is relaxed. All one needs is the assumption that nonzero coefficients of $\beta$ are "sufficiently" large. One possibility is hard-thresholding of the obtained coefficients, neglecting variables with very small coefficients. This effect has already been observed empirically in [33]. Other possibilities include soft-thresholding and relaxation methods such as the Gauss–Dantzig selector [5], the relaxed Lasso [22] with an additional thresholding step or the adaptive Lasso of Zou [42].

DEFINITION 5 (Hard-thresholded Lasso estimator). Let, for each $x \in \mathbb{R}^p$, the quantity $1\{|x| \geq c\}$ be a $p_n$-dimensional vector which is, componentwise, equal to 1 if $|x_k| \geq c$ and 0 otherwise. For a given sequence $t_n$, the hard-thresholded Lasso estimator $\hat{\beta}^{ht,\lambda}$ is defined as

$$\hat{\beta}^{ht,\lambda} = \hat{\beta}^\lambda 1\{\hat{\beta}^\lambda \geq \sigma t_n \sqrt{\log p_n/n}\}.$$

The sequence $t_n$ can be chosen freely. We start with a corollary that follows directly from Theorem 1, stating that the hard-thresholded Lasso estimator (unlike the un-thresholded estimator) is *sign consistent* under regularity assumptions that are weaker than the *irrepresentable condition* needed for sign-consistency of the ordinary Lasso estimator.

COROLLARY 2 (Sign consistency by hard thresholding). *Assume the incoherent design assumption ([6]) holds and the sparsity of $\beta$ fulfills $s_n = o(t_n^2 e_n^{-4})$ for $n \to \infty$. Assume furthermore*

$$\min_{k\,:\,\beta_k \neq 0} |\beta_k| \gg \sigma t_n \sqrt{\log p_n/n}, \qquad n \to \infty.$$

*Under a choice $\lambda \propto \sigma e_n \sqrt{n \log p_n}$, the hard-thresholded Lasso estimator of Definition [5] is then sign-consistent and*

$$P\{\text{sign}(\hat{\beta}^{ht,\lambda}) = \text{sign}(\beta)\} \to 1 \qquad as\ n \to \infty.$$

The proof follows from the results of Theorem [1]. The bound ([12]) on the $\ell_2$-distance, derived from Theorem [1], gives then trivially the identical bound on the squared $\ell_\infty$-distance between $\hat{\beta}^\lambda$ and $\beta$. The result follows by observing that $1/\phi_{\max} = O(1)$ and the fact that $\ell_\infty$ error is a smaller order of the lower bound on the size of nonzero $\beta$'s due to assumptions of incoherent design and $s_n = o(t_n^2 e_n^{-4})$. When choosing a suitable value of the cut-off parameter $t_n$, one is faced with a trade-off. Choosing larger values of the cut-off $t_n$ places a stricter condition on the minimal nonzero value of $\beta$, while smaller values of $t_n$ relax this assumption, yet require the vector $\beta$ to be sparser.

The result mainly implies that *sign-consistency* can be achieved with the hard-thresholded Lasso estimator under much weaker consistency requirements than with the ordinary Lasso estimator. As discussed previously, the ordinary Lasso estimator is only sign consistent if the *irrepresentable condition* or, equivalently, *neighborhood stability* is fulfilled [23, 41, 42]. This is a considerably stronger assumption than the incoherence assumption above. In either case, a similar assumption on the rate of decay of the minimal nonzero components is needed.

In conclusion, even though one cannot achieve *sign consistency* in general with just a single Lasso estimation, it can be achieved in a two-stage procedure.

**3. Proof of Theorem 1.** Let $\beta^\lambda$ be the estimator under the absence of noise, that is, $\beta^\lambda = \hat{\beta}^{\lambda,0}$, where $\hat{\beta}^{\lambda,\xi}$ is defined as in (15). The $\ell_2$-distance can then be bounded by $\|\hat{\beta}^\lambda - \beta\|_{\ell_2}^2 \leq 2\|\hat{\beta}^\lambda - \beta^\lambda\|_{\ell_2}^2 + 2\|\beta^\lambda - \beta\|_{\ell_2}^2$. The first term on the right-hand side represents the variance of the estimation, while the second term represents the bias. The bias contribution follows directly from Lemma 2 below. The bound on the variance term follows by Lemma 6 below.

*De-noised response.* Before starting, it is useful to define a de-noised response. Define for $0 < \xi < 1$ the de-noised version of the response variable,

$$Y(\xi) = X\beta + \xi\varepsilon. \tag{14}$$

We can regulate the amount of noise with the parameter $\xi$. For $\xi = 0$, only the signal is retained. The original observations with the full amount of noise are recovered for $\xi = 1$. Now consider for $0 \leq \xi \leq 1$ the estimator $\hat{\beta}^{\lambda,\xi}$,

$$\hat{\beta}^{\lambda,\xi} = \arg\min_\beta \|Y(\xi) - X\beta\|_{\ell_2}^2 + \lambda\|\beta\|_{\ell_1}. \tag{15}$$

The ordinary Lasso estimate is recovered under the full amount of noise so that $\hat{\beta}^{\lambda,1} = \hat{\beta}^\lambda$. Using the notation from the previous results, we can write for the estimate in the absence of noise, $\hat{\beta}^{\lambda,0} = \beta^\lambda$. The definition of the de-noised version of the Lasso estimator will be helpful for the proof as it allows to characterize the variance of the estimator.

3.1. *Part* I *of proof: bias.* Let $K$ be the set of nonzero elements of $\beta$, that is, $K = \{k : \beta_k \neq 0\}$. The cardinality of $K$ is again denoted by $s = s_n$. For the following, let $\beta^\lambda$ be the estimator $\hat{\beta}^{\lambda,0}$ under the absence of noise, as defined in (15). The solution $\beta^\lambda$ can, for each value of $\lambda$, be written as $\beta^\lambda = \beta + \gamma^\lambda$, where

$$\gamma^\lambda = \arg\min_{\zeta \in \mathbb{R}^p} f(\zeta). \tag{16}$$

The function $f(\zeta)$ is given by

$$f(\zeta) = n\zeta^T C \zeta + \lambda \sum_{k \in K^c} |\zeta_k| + \lambda \sum_{k \in K} (|\beta_k + \zeta_k| - |\beta_k|). \tag{17}$$

The vector $\gamma^\lambda$ is the bias of the Lasso estimator. We derive first a bound on the $\ell_2$-norm of $\gamma^\lambda$.

LEMMA 1. *Assume incoherent design as in (6) with a sparsity multiplier sequence $e_n$. The $\ell_2$-norm of $\gamma^\lambda$, as defined in (16), is then bounded for sufficiently large values of $n$ by*

$$\|\gamma^\lambda\|_{\ell_2} \leq 17.5 \frac{\lambda}{n} \frac{\sqrt{s_n}}{\phi_{\min}(e_n s_n)}. \tag{18}$$

PROOF. We write in the following $\gamma$ instead of $\gamma^\lambda$ for notational simplicity. Let $\gamma(K)$ be the vector with coefficients $\gamma_k(K) = \gamma_k 1\{k \in K\}$, that is, $\gamma(K)$ is the bias of the truly nonzero coefficients. Analogously, let $\gamma(K^c)$ be the bias of the truly zero coefficients with $\gamma_k(K^c) = \gamma_k 1\{k \notin K\}$. Clearly, $\gamma = \gamma(K) + \gamma(K^c)$. The value of the function $f(\zeta)$, as defined in (17), is 0 if setting $\zeta = 0$. For the true solution $\gamma^\lambda$, it follows hence that $f(\gamma^\lambda) \le 0$. Hence, using that $\zeta^T C \zeta \ge 0$ for any $\zeta$,

$$(19) \qquad \|\gamma(K^c)\|_{\ell_1} = \sum_{k \in K^c} |\zeta_k| \le \left| \sum_{k \in K} (|\beta_k + \zeta_k| - |\beta_k|) \right| \le \|\gamma(K)\|_{\ell_1}.$$

As $\|\gamma(K)\|_{\ell_0} \le s_n$, it follows that $\|\gamma(K)\|_{\ell_1} \le \sqrt{s_n}\|\gamma(K)\|_{\ell_2} \le \sqrt{s_n}\|\gamma\|_{\ell_2}$ and hence, using (19),

$$(20) \qquad \|\gamma\|_{\ell_1} \le 2\sqrt{s_n}\|\gamma\|_{\ell_2}.$$

This result will be used further below. We use now again that $f(\gamma^\lambda) \le 0$ [as $\zeta = 0$ yields the upper bound $f(\zeta) = 0$]. Using the previous result that $\|\gamma(K)\|_{\ell_1} \le \sqrt{s_n}\|\gamma\|_{\ell_2}$, and ignoring the nonnegative term $\|\gamma(K^c)\|_{\ell_1}$, it follows that

$$(21) \qquad n\gamma^T C \gamma \le \lambda \sqrt{s_n}\|\gamma\|_{\ell_2}.$$

Consider now the term $\gamma^T C \gamma$. Bounding this term from below and plugging the result into (21) will yield the desired upper bound on the $\ell_2$-norm of $\gamma$. Let $|\gamma_{(1)}| \ge |\gamma_{(2)}| \ge \cdots \ge |\gamma_{(p)}|$ be the ordered entries of $\gamma$.

Let $u_n$ for $n \in \mathbb{N}$ be a sequence of positive integers, to be chosen later, and define the set of the "$u_n$-largest coefficients" as $U = \{k : |\gamma_k| \ge |\gamma_{(u_n)}|\}$. Define analogously to above the vectors $\gamma(U)$ and $\gamma(U^c)$ by $\gamma_k(U) = \gamma_k 1\{k \in U\}$ and $\gamma_k(U^c) = \gamma_k 1\{k \notin U\}$. The quantity $\gamma^T C \gamma$ can be written as $\gamma^T C \gamma = \|a + b\|_{\ell_2}^2$, where $a := n^{-1/2} X \gamma(U)$ and $b := n^{-1/2} X \gamma(U^c)$. Then

$$(22) \qquad \gamma^T C \gamma = \|a + b\|_{\ell_2}^2 \ge (\|a\|_{\ell_2} - \|b\|_{\ell_2})^2.$$

Before proceeding, we need to bound the norm $\|\gamma(U^c)\|_{\ell_2}$ as a function of $u_n$. Assume for the moment that the $\ell_1$-norm $\|\gamma\|_{\ell_1}$ is identical to some $\ell > 0$. Then it holds for every $k = 1, \ldots, p$ that $\gamma_{(k)} \le \ell/k$. Hence,

$$(23) \qquad \|\gamma(U^c)\|_{\ell_2}^2 \le \|\gamma\|_{\ell_1}^2 \sum_{k=u_n+1}^{p} \frac{1}{k^2} \le (4 s_n \|\gamma\|_{\ell_2}^2) \frac{1}{u_n},$$

having used the result (20) from above that $\|\gamma\|_{\ell_1} \le 2\sqrt{s_n}\|\gamma\|_{\ell_2}$. As $\gamma(U)$ has by definition only $u_n$ nonzero coefficients,

$$(24) \qquad \begin{aligned} \|a\|_{\ell_2}^2 &= \|\gamma(U)^T C \gamma(U)\|_{\ell_2}^2 \ge \phi_{\min}(u_n) \|\gamma(U)\|_{\ell_2}^2 \\ &\ge \phi_{\min}(u_n) \left(1 - \frac{4 s_n}{u_n}\right) \|\gamma\|_{\ell_2}^2, \end{aligned}$$

having used (23) and $\|\gamma(U)\|_{\ell_2}^2 = \|\gamma\|_{\ell_2}^2 - \|\gamma(U^c)\|_{\ell_2}^2$. As $\gamma(U^c)$ has at most $\min\{n,p\}$ nonzero coefficients and using again (23),

$$(25) \quad \|b\|_{\ell_2}^2 = \|\gamma(U^c)^T C \gamma(U^c)\|_{\ell_2}^2 \leq \phi_{\max} \|\gamma(U^c)\|_{\ell_2}^2 \leq \phi_{\max} \frac{4s_n}{u_n} \|\gamma\|_{\ell_2}^2.$$

Using (24) and (25) in (22), together with $\phi_{\max} \geq \phi_{\min}(u_n)$,

$$(26) \quad \gamma^T C \gamma \geq \phi_{\min}(u_n) \|\gamma\|_{\ell_2}^2 \left(1 - 4\sqrt{\frac{s_n \phi_{\max}}{u_n \phi_{\min}(u_n)}}\right).$$

Choosing for $u_n$ the *sparsity multiplier* sequence, as defined in (6), times the sparsity $s_n$, so that $u_n = e_n s_n$ it holds that $s_n \phi_{\max}/(e_n s_n \phi_{\min}(e_n s_n)) < 1/18$ and hence also that $s_n \phi_{\max}/(e_n s_n \phi_{\min}(e_n^2 s_n)) < 1/18$, since $\phi_{\min}(e_n^2 s_n) \leq \phi_{\min}(e_n s_n)$. Thus the right-hand side in (26) is bounded from below by $18\phi_{\min}(e_n s_n)\|\gamma\|_{\ell_2}^2$ since $(1 - 4/\sqrt{18}) \leq 17.5$. Using the last result together with (21), which says that $\gamma^T C \gamma \leq n^{-1} \lambda \sqrt{s_n} \|\gamma\|_{\ell_2}$, it follows that for large $n$,

$$\|\gamma\|_{\ell_2} \leq 17.5 \frac{\lambda}{n} \frac{\sqrt{s_n}}{\phi_{\min}(e_n s_n)},$$

which completes the proof. $\square$

LEMMA 2. *Under the assumptions of Theorem 1, the bias $\|\gamma^\lambda\|_{\ell_2}^2$ is bounded by*

$$\|\gamma^\lambda\|_{\ell_2}^2 \leq (17.5)^2 \sigma^2 \frac{s_n \log p_n}{n} \frac{e_n^2}{\phi_{\min}^2(e_n^2 s_n)}.$$

PROOF. This is an immediate consequence of Lemma 1. Plugging the penalty sequence $\lambda \propto \sigma \sqrt{n \log p_n} e_n$ into (18), the results follows by the inequality $\phi_{\min}(e_n s_n) \geq \phi_{\min}(e_n^2 s_n)$, having used that, by its definition in (6), $e_n$ is necessarily larger than 1. $\square$

3.2. *Part* II *of proof: variance.* The proof for the variance part needs two steps. First, a bound on the variance is derived, which is a function of the number of active variables. In a second step, the number of active variables will be bounded, taking into account also the bound on the bias derived above.

*Variance of restricted OLS.* Before considering the Lasso estimator, a trivial bound is shown for the variance of a restricted OLS estimation. Let $\hat{\theta}^M \in \mathbb{R}^p$ be, for every subset $M \subseteq \{1,\ldots,p\}$ with $|M| \leq n$, the restricted OLS-estimator of the noise vector $\varepsilon$,

$$(27) \qquad \hat{\theta}^M = (X_M^T X_M)^{-1} X_M^T \varepsilon.$$

First, we bound the $\ell_2$-norm of this estimator. The result is useful for bounding the variance of the final estimator, based on the derived bound on the number of active variables.

LEMMA 3. *Let $\overline{m}_n$ be a sequence with $\overline{m}_n = o(n)$ and $\overline{m}_n \to \infty$ for $n \to \infty$. If $p_n \to \infty$, it holds with probability converging to 1 for $n \to \infty$*

$$\max_{M:|M|\leq \overline{m}_n} \|\hat{\theta}^M\|_{\ell_2}^2 \leq \frac{2\log p_n}{n} \frac{\overline{m}_n}{\phi_{\min}^2(\overline{m}_n)} \sigma^2.$$

The $\ell_2$-norm of the restricted estimator $\hat{\theta}^M$ is thus bounded uniformly over all sets $M$ with $|M| \leq \overline{m}_n$.

PROOF OF LEMMA 3. It follows directly from the definition of $\hat{\theta}^M$ that, for every $M$ with $|M| \leq \overline{m}_n$,

(28) $$\|\hat{\theta}^M\|_{\ell_2}^2 \leq \frac{1}{n^2 \phi_{\min}^2(\overline{m}_n)} \|X_M^T \varepsilon\|_{\ell_2}^2.$$

It remains to be shown that, for $n \to \infty$, with probability converging to 1,

$$\max_{M:|M|\leq \overline{m}_n} \|X_M^T \varepsilon\|_{\ell_2}^2 \leq 2\log p_n \sigma^2 \overline{m}_n n.$$

As $\varepsilon_i \sim \mathcal{N}(0, \sigma^2)$ for all $i = 1, \ldots, n$, it holds with probability converging to 1 for $n \to \infty$, by Bonferroni's inequality that $\max_{k \leq p_n} |X_k^T \varepsilon|^2$ is bounded from above by $2\log p_n \sigma^2 n$. Hence, with probability converging to 1 for $n \to \infty$,

(29) $$\max_{M:|M|\leq \overline{m}_n} \|X_M^T \varepsilon\|_{\ell_2}^2 \leq \overline{m}_n \max_{k \leq p_n} |X_k^T \varepsilon|^2 \leq 2\log p_n \sigma^2 n \overline{m}_n,$$

which completes the proof. □

*Variance of estimate is bounded by restricted OLS variance.* We show that the variance of the Lasso estimator can be bounded by the variances of restricted OLS estimators, using bounds on the number of active variables.

LEMMA 4. *If, for a fixed value of $\lambda$, the number of active variables of the de-noised estimators $\hat{\beta}^{\lambda,\xi}$ is for every $0 \leq \xi \leq 1$ bounded by $m$, then*

(30) $$\sup_{0 \leq \xi \leq 1} \|\hat{\beta}^{\lambda,0} - \hat{\beta}^{\lambda,\xi}\|_{\ell_2}^2 \leq \max_{M:|M|\leq m} \|\hat{\theta}^M\|_{\ell_2}^2.$$

PROOF. The key in the proof is that the solution path of $\hat{\beta}^{\lambda,\xi}$, if increasing the value of $\xi$ from 0 to 1, can be expressed piecewise in terms of the restricted OLS solution. It will be obvious from the proof that it is sufficient to show the claim for $\xi = 1$ in the term on the r.h.s. of (30).

The set $M(\xi)$ of active variables is the set with maximal absolute gradient,

$$M(\xi) = \{k : |G_k^{\lambda,\xi}| = \lambda\}.$$

Note that the estimator $\hat{\beta}^{\lambda,\xi}$ and also the gradient $G_k^{\lambda,\xi}$ are continuous functions in both $\lambda$ and $\xi$ [11]. Let $0 = \xi_1 < \xi_2 < \cdots < \xi_{L+1} = 1$ be the points of discontinuity of $M(\xi)$. At these locations, variables either join the active set or are dropped from the active set.

Fix some $j$ with $1 \leq j \leq J$. Denote by $M_j$ the set of active variables $M(\xi)$ for any $\xi \in (\xi_j, \xi_{j+1})$. We show in the following that the solution $\hat{\beta}^{\lambda,\xi}$ is for all $\xi$ in the interval $(\xi_j, \xi_{j+1})$ given by

$$(31) \qquad \forall \xi \in (\xi_j, \xi_{j+1}) : \hat{\beta}^{\lambda,\xi} = \hat{\beta}^{\lambda,\xi_j} + (\xi - \xi_j)\hat{\theta}^{M_j},$$

where $\hat{\theta}^{M_j}$ is the restricted OLS estimator of noise, as defined in (27). The local effect of increased noise (larger value of $\xi$) on the estimator is thus to shift the coefficients of the active set of variables along the least squares direction.

Once (31) is shown, the claim follows by piecing together the piecewise linear parts and using continuity of the solution as a function of $\xi$ to obtain

$$\|\hat{\beta}^{\lambda,0} - \hat{\beta}^{\lambda,1}\|_{\ell_2} \leq \sum_{j=1}^{J} \|\hat{\beta}^{\lambda,\xi_j} - \hat{\beta}^{\lambda,\xi_{j+1}}\|_{\ell_2}$$

$$\leq \max_{M:|M|\leq m} \|\hat{\theta}^M\|_{\ell_2} \sum_{j=1}^{J}(\xi_{j+1} - \xi_j) = \max_{M:|M|\leq m} \|\hat{\theta}^M\|_{\ell_2}.$$

It thus remains to show (31). A necessary and sufficient condition for $\hat{\beta}^{\lambda,\xi}$ with $\xi \in (\xi_j, \xi_{j+1})$ to be a valid solution is that for all $k \in M_j$ with nonzero coefficient $\hat{\beta}_k^{\lambda,\xi} \neq 0$, the gradient is equal to $\lambda$ times the negative sign,

$$(32) \qquad G_k^{\lambda,\xi} = -\lambda \operatorname{sign}(\hat{\beta}_k^{\lambda,\xi}),$$

that for all variables with $k \in M_j$ with zero coefficient $\hat{\beta}_k^{\lambda,\xi} = 0$ the gradient is equal in absolute value to $\lambda$

$$(33) \qquad |G_k^{\lambda,\xi}| = \lambda$$

and for variables $k \notin M_j$ not in the active set,

$$(34) \qquad |G_k^{\lambda,\xi}| < \lambda.$$

These conditions are a consequence of the requirement that the subgradient of the loss function contains 0 for a valid solution.

Note that the gradient of the active variables in $M_j$ is unchanged if replacing $\xi \in (\xi_j, \xi_{j+1})$ by some $\xi' \in (\xi_j, \xi_{j+1})$ and replacing $\hat{\beta}^{\lambda,\xi}$ by $\hat{\beta}^{\lambda,\xi} + (\xi' - \xi)\hat{\theta}^{M_j}$. That is, for all $k \in M_j$,

$$(Y(\xi) - X\hat{\beta}^{\lambda,\xi})^T X_k = \{Y(\xi') - X(\hat{\beta}^{\lambda,\xi} + (\xi' - \xi)\hat{\theta}^{M_j})\}^T X_k,$$

as the difference of both sides is equal to $(\xi' - \xi)\{(\varepsilon - X\hat{\theta}^{M_j})^T X_k\}$, and $(\varepsilon - X\hat{\theta}^{M_j})^T X_k = 0$ for all $k \in M_j$, as $\hat{\theta}^{M_j}$ is the OLS of $\varepsilon$, regressed on the variables in $M_j$. Equalities (32) and (33) are thus fulfilled for the solution and it remains to show that (34) also holds. For sufficiently small values of $\xi' - \xi$, inequality (34) is clearly fulfilled for continuity reasons. Note that if $|\xi' - \xi|$ is large enough such that for one variable $k \notin M_j$ inequality (34) becomes an equality, then the set of active variables changes and thus either $\xi' = \xi_{j+1}$ or $\xi' = \xi_j$. We have thus shown that the solution $\hat{\beta}^{\lambda,\xi}$ can for all $\xi \in (\xi_j, \xi_{j+1})$ be written as

$$\hat{\beta}^{\lambda,\xi} = \hat{\beta}^{\lambda,\xi_j} + (\xi - \xi_j)\hat{\theta}^{M_j},$$

which proves (31) and thus completes the proof. □

*A bound on the number of active variables.* A decisive part in the variance of the estimator is determined by the number of selected variables. Instead of directly bounding the number of selected variables, we derive bounds for the number of *active variables*. As any variable with a nonzero regression coefficient is also an *active variable*, these bounds lead trivially to bounds for the number of selected variables.

Let $\mathcal{A}_\lambda$ be the set of *active variables*,

$$\mathcal{A}_\lambda = \{k : |G_k^\lambda| = \lambda\}.$$

Let $\mathcal{A}_{\lambda,\xi}$ be the set of *active variables* of the de-noised estimator $\hat{\beta}^{\lambda,\xi}$, as defined in (15). The number of selected variables (variables with a nonzero coefficient) is at most as large as the number of *active variables*, as any variable with a nonzero estimated coefficient has to be an *active variable* [25].

LEMMA 5. *For $\lambda \geq \sigma e_n \sqrt{n \log p_n}$, the maximal number $\sup_{0 \leq \xi \leq 1} |\mathcal{A}_{\lambda,\xi}|$ of active variables is bounded, with probability converging to 1 for $n \to \infty$, by*

$$\sup_{0 \leq \xi \leq 1} |\mathcal{A}_{\lambda,\xi}| \leq e_n^2 s_n.$$

PROOF. Let $R^{\lambda,\xi}$ be the residuals of the de-noised estimator (15), $R^{\lambda,\xi} = Y - X\hat{\beta}^{\lambda,\xi}$. For any $k$ in the $|\mathcal{A}_{\lambda,\xi}|$-dimensional space spanned by the *active variables*,

(35) $$|X_k^T R^{\lambda,\xi}| = \lambda.$$

Adding up, it follows that for all $0 \leq \xi \leq 1$,

$$|\mathcal{A}_{\lambda,\xi}|\lambda^2 = \|X^T_{\mathcal{A}_{\lambda,\xi}} R^{\lambda,\xi}\|^2_{\ell_2}. \tag{36}$$

The residuals can for all values $0 \leq \xi \leq 1$ be written as the sum of two terms, $R^{\lambda,\xi} = X(\beta - \hat{\beta}^{\lambda,\xi}) + \xi\varepsilon$. Equality (36) can now be transformed into the inequality,

$$|\mathcal{A}_{\lambda,\xi}|\lambda^2 \leq (\|X^T_{\mathcal{A}_{\lambda,\xi}} X(\beta - \hat{\beta}^{\lambda,\xi})\|_{\ell_2} + \xi^2\|X^T_{\mathcal{A}_{\lambda,\xi}}\varepsilon\|_{\ell_2})^2 \tag{37}$$

$$\leq (\|X^T_{\mathcal{A}_{\lambda,\xi}} X(\beta - \hat{\beta}^{\lambda,\xi})\|_{\ell_2} + \|X^T_{\mathcal{A}_{\lambda,\xi}}\varepsilon\|_{\ell_2})^2. \tag{38}$$

Denote by $\tilde{m}$ the supremum of $|\mathcal{A}_{\lambda,\xi}|$ over all values of $0 \leq \xi \leq 1$. Using the same argument as in the derivation of (29), the term $\sup_{0 \leq \xi \leq 1} \|X^T_{\mathcal{A}_{\lambda,\xi}}\varepsilon\|^2_{\ell_2}$ is of order $o_p(\tilde{m} n \log p_n)$ as long as $p_n \to \infty$ for $n \to \infty$. For sufficiently large $n$ it holds thus, using $\lambda \geq \sigma e_n \sqrt{n \log p_n}$, that $\sup_{0 \leq \xi \leq 1} \|X^T_{\mathcal{A}_{\lambda,\xi}}\varepsilon\|_{\ell_2}/(\tilde{m}\lambda^2)^{1/2} \leq \eta$ for any $\eta > 0$. Dividing by $\lambda^2$, (37) implies then, with probability converging to 1,

$$\tilde{m} \leq \sup_{0 \leq \xi \leq 1} (\lambda^{-1}\|X^T_{\mathcal{A}_{\lambda,\xi}} X(\beta - \hat{\beta}^{\lambda,\xi})\|_{\ell_2} + \eta\sqrt{\tilde{m}})^2. \tag{39}$$

Now turning to the right-hand side, it trivially holds for any value of $0 \leq \xi \leq 1$ that $|\mathcal{A}_{\lambda,\xi}| \leq \min\{n,p\}$. On the other hand, $X(\beta - \hat{\beta}^{\lambda,\xi}) = X_{\mathcal{B}_{\lambda,\xi}}(\beta - \hat{\beta}^{\lambda,\xi})$, where $\mathcal{B}_{\lambda,\xi} := \mathcal{A}_{\lambda,\xi} \cup \{k : \beta_k \neq 0\}$, as the difference vector $\beta - \hat{\beta}^{\lambda,\xi}$ has nonzero entries only in the set $\mathcal{B}_{\lambda,\xi}$. Thus

$$\|X^T_{\mathcal{A}_{\lambda,\xi}} X(\beta - \hat{\beta}^{\lambda,\xi})\|^2_{\ell_2} \leq \|X^T_{\mathcal{B}_{\lambda,\xi}} X_{\mathcal{B}_{\lambda,\xi}}(\beta - \hat{\beta}^{\lambda,\xi})\|^2_{\ell_2}.$$

Using additionally $|\mathcal{B}_{\lambda,\xi}| \leq s_n + \min\{n,p\}$, it follows that

$$\|X^T_{\mathcal{A}_{\lambda,\xi}} X(\beta - \hat{\beta}^{\lambda,\xi})\|^2_{\ell_2} \leq n^2 \phi^2_{\max}\|(\beta - \hat{\beta}^{\lambda,\xi})\|^2_{\ell_2}.$$

Splitting the difference $\beta - \hat{\beta}^{\lambda,\xi}$ into $(\beta - \beta^{\lambda}) + (\beta^{\lambda} - \hat{\beta}^{\lambda,\xi})$, where $\beta^{\lambda} = \hat{\beta}^{\lambda,0}$ is again the population version of the Lasso estimator, it holds for any $\eta > 0$, using (39), that with probability converging to 1 for $n \to \infty$,

$$\tilde{m} \leq \left(n\lambda^{-1}\phi_{\max}\|\beta - \beta^{\lambda}\|_{\ell_2} + n\lambda^{-1}\phi_{\max} \sup_{0 \leq \xi \leq 1}\|\hat{\beta}^{\lambda,0} - \hat{\beta}^{\lambda,\xi}\|_{\ell_2} + \eta\sqrt{\tilde{m}}\right)^2. \tag{40}$$

Using Lemmas 3 and 4, the variance term $n^2\phi^2_{\max}\sup_{0 \leq \xi \leq 1}\|\hat{\beta}^{\lambda,0} - \hat{\beta}^{\lambda,\xi}\|^2_{\ell_2}$ is bounded by $o_p\{n\tilde{m}\log p_n \phi^2_{\max}/\phi^2_{\min}(\tilde{m})\}$. Define, implicitly, a sequence $\tilde{\lambda} = \sigma\sqrt{n \log p_n}(\phi_{\max}/\phi_{\min}(\tilde{m}))$. For any sequence $\lambda$ with $\liminf_{n\to\infty} \lambda/\tilde{\lambda} > 0$, the term $n^2\lambda^{-2}\phi^2_{\max}\sup_{0 \leq \xi \leq 1}\|\hat{\beta}^{\lambda,0} - \hat{\beta}^{\lambda,\xi}\|^2_{\ell_2}$ is then of order $o_p(\tilde{m})$. Using furthermore the bound on the bias from Lemma 1, it holds with probability

converging to 1, for $n \to \infty$ for any sequence $\lambda$ with $\liminf_{n\to\infty} \lambda/\tilde{\lambda} > 0$ and any $\eta > 0$ that

$$\tilde{m} \leq (n\lambda^{-1}\phi_{\max}\|\beta - \beta^\lambda\|_{\ell_2} + 2\eta\sqrt{\tilde{m}})^2 \leq \left(17.5\phi_{\max}\frac{\sqrt{s_n}}{\phi_{\min}(e_n s_n)} + 2\eta\sqrt{\tilde{m}}\right)^2.$$

Choosing $\eta = 0.013$ implies, for an inequality of the form $a^2 \leq (x + 2\eta a)^2$, that $a \leq (18/17.5)x$. Hence, choosing this value of $\eta$, it follows from the equation above that, with probability converging to 1 for $n \to \infty$,

$$\tilde{m} \leq 18^2 \phi_{\max}^2 \frac{s_n}{\phi_{\min}^2(e_n s_n)} = e_n^2 s_n \left(\frac{18\phi_{\max}}{e_n \phi_{\min}(e_n s_n)}\right)^2 \leq e_n^2 s_n,$$

having used the definition of the *sparsity multiplier* in (6). We can now see that the requirement on $\lambda$, namely $\liminf_{n\to\infty} \lambda/\tilde{\lambda} > 0$, is fulfilled if $\lambda \geq \sigma e_n \sqrt{n \log p_n}$, which completes the proof. □

Finally, we use Lemmas 3, 4 and 5 to show the bound on the variance of the estimator.

LEMMA 6. *Under the conditions of Theorem 1, with probability converging to 1 for $n \to \infty$,*

$$\|\beta^\lambda - \hat{\beta}^\lambda\|_{\ell_2}^2 \leq 2\sigma^2 \frac{s_n \log p_n}{n} \frac{e_n^2}{\phi_{\min}^2(e_n^2 s_n)}.$$

The proof follows immediately from Lemmas 3 and 4 when inserting the bound on the number of active variables obtained in Lemma 5.

**4. Numerical illustration: frequency detection.** Instead of extensive numerical simulations, we would like to illustrate a few aspects of Lasso-type variable selection if the *irrepresentable condition* is not fulfilled. We are not making claims that the Lasso is superior to other methods for high-dimensional data. We merely want to draw attention to the fact that (a) the Lasso might not be able to select the correct variables but (b) comes nevertheless close to the true vector in an $\ell_2$-sense.

An illustrative example is frequency detection. It is of interest in some areas of the physical sciences to accurately detect and resolve frequency components; two examples are variable stars [27] and detection of gravitational waves [6, 32]. A nonparametric approach is often most suitable for fitting of the involved periodic functions [15]. However, we assume here for simplicity that the observations $Y = (Y_1, \ldots, Y_n)$ at time points $t = (t_1, \ldots, t_n)$ are of the form

$$Y_i = \sum_{\omega \in \Omega} \beta_\omega \sin(2\pi\omega t_i + \phi_\omega) + \varepsilon_i,$$

where $\Omega$ contains the set of fundamental frequencies involved, and $\varepsilon_i$ for $i = 1, \ldots, n$ is independently and identically distributed noise with $\varepsilon_i \sim \mathcal{N}(0, \sigma^2)$. To simplify the problem even more, we assume that the phases are known to be zero, $\phi_\omega = 0$ for all $\omega \in \Omega$. Otherwise one might like to employ the Group Lasso [37], grouping together the sine and cosine part of identical frequencies.

It is of interest to resolve closely adjacent spectral lines [16] and we will work in this setting in the following. We choose for the experiment $n = 200$ evenly spaced observation times. There are supposed to be two closely adjacent frequencies with $\omega_1 = 0.0545$ and $\omega_2 = 0.0555 = \omega_1 + 1/300$, both entering with $\beta_{\omega_1} = \beta_{\omega_2} = 1$. As we have the information that the phase is zero for all frequencies, the predictor variables are given by all sine-functions with frequencies evenly spaced between $1/200$ and $1/2$, with a spacing of $1/600$ between adjacent frequencies.

In the chosen setting, the *irrepresentable condition* is violated for the frequency $\omega_m = (\omega_1 + \omega_2)/2$. Even in the absence of noise, this *resonance frequency* is included in the Lasso-estimate for all positive penalty parameters, as can be seen from the results further below. As a consequence of a violated *irrepresentable condition*, the largest peak in the periodogram is in general obtained for the *resonance frequency*. In Figure 1 we show the periodogram [28] under a moderate noise level $\sigma = 0.2$. The periodogram shows the amount of energy in each frequency, and is defined through the function

$$\Delta E(\omega) = \sum_i Y_i^2 - \sum_i (Y_i - \hat{Y}_i^{(\omega)})^2,$$

where $\hat{Y}^{(\omega)}$ is the least squares fit of the observations $Y$, using only sine and cosine functions with frequency $\omega$ as two predictor variables. There is clearly a peak at frequency $\omega_m$. As can be seen in the close-up around $\omega_m$, it is not immediately obvious from the periodogram that there are *two* frequencies at frequencies $\omega_1$ and $\omega_2$. As said above, the *irrepresentable condition* is violated for the *resonance frequency* and it is of interest to see which frequencies are picked up by the Lasso estimator.

The results are shown in Figures 2 and 3. Figure 3 highlights that the two true frequencies are with high probability picked up by the Lasso. The *resonance frequency* is also selected with high probability, no matter how the penalty is chosen. This result could be expected as the *irrepresentable condition* is violated and the estimator can thus not be *sign consistent*. We expect from the theoretical results in this manuscript that the coefficient of the falsely selected *resonance frequency* is very small if the penalty parameter is chosen correctly. And it can indeed be seen in Figure 2 that the coefficients of the true frequencies are much larger than the coefficient of the *resonance frequency* for an appropriate choice of the penalty parameter.

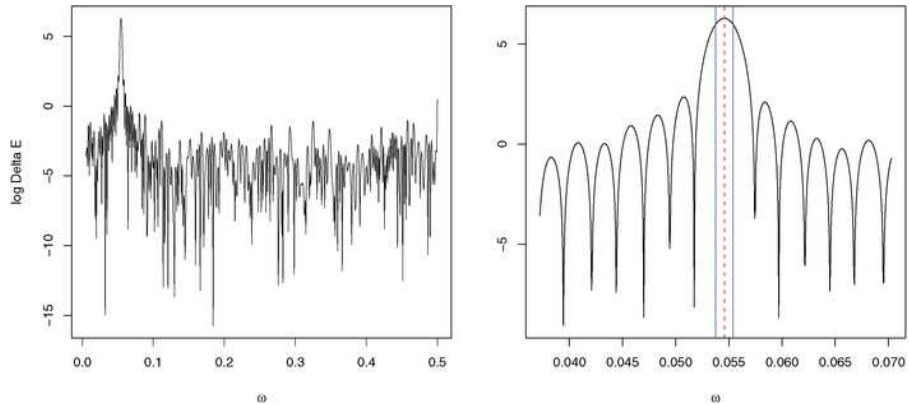

Fig. 1. *The energy $\log \Delta E(\omega)$ for a noise level $\sigma = 0.2$ is shown on the left for a range of frequencies $\omega$. A close-up of the region around the peak is shown on the right. The two frequencies $\omega_1$ and $\omega_2$ are marked with solid vertical lines, while the resonance frequency $(\omega_1 + \omega_2)/2$ is shown with a broken vertical line.*

These results reinforce our conclusion that the Lasso might not be able to pick up the correct sparsity pattern, but delivers nevertheless useful approximations as falsely selected variables are chosen only with a very small coefficient; this behavior is typical and expected from the results of Theorem 1. Falsely selected coefficients can thus be removed in a second step, either by thresholding variables with small coefficients or using other relaxation techniques. In any case, it is reassuring to know that all important variables are included in the Lasso estimate.

**5. Concluding remarks.** It has recently been discovered that the Lasso cannot recover the correct sparsity pattern in certain circumstances, even not asymptotically for $p_n$ fixed and $n \to \infty$. This sheds a little doubt on whether the Lasso is a good method for identification of sparse models for both low- and high-dimensional data.

Here we have shown that the Lasso can continue to deliver good approximations to sparse coefficient vectors $\beta$ in the sense that the $\ell_2$-difference $\|\beta - \hat{\beta}^\lambda\|_{\ell_2}$ vanishes for large sample sizes $n$, even if it fails to discover the correct sparsity pattern. The conditions needed for a good approximation in the $\ell_2$-sense are weaker than the *irrepresentable condition* needed for *sign consistency*. We pointed out that the correct sparsity pattern could be recovered in a two-stage procedure when the true coefficients are not too small. The first step consists in a regular Lasso fit. Variables with small absolute coefficients are then removed from the model in a second step.

We derived possible scenarios under which $\ell_2$-consistency in the sense of (4) can be achieved as a function of the sparsity of the vector $\beta$, the number of samples and the number of variables. Under the condition that sparse minimal eigenvalues are not decaying too fast in some sense, the requirement for $\ell_2$-consistency is (ignoring $\log n$ factors)

$$\frac{s_n \log p_n}{n} \to 0 \qquad \text{as } n \to \infty.$$

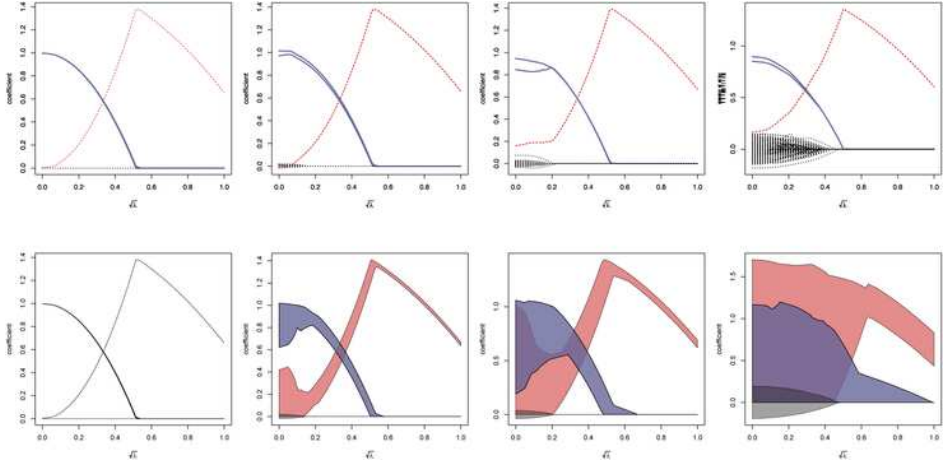

FIG. 2. *An example where the Lasso is bound to select wrong variables, while being a good approximation to the true vector in the $\ell_2$-sense. Top row: The noise level increases from left to right as $\sigma = 0, 0.1, 0.2, 1$. For one run of the simulation, paths of the estimated coefficients are shown as a function of the square root $\sqrt{\lambda}$ of the penalty parameter. The actually present signal frequencies $\omega_1$ and $\omega_2$ are shown as solid lines, the resonance frequency as a broken line, and all other frequencies are shown as dotted lines. Bottom row: The shaded areas contain, for 90% of all simulations, the regularization paths of the signal frequencies (region with solid borders), resonance frequency (area with broken borders) and all other frequencies (area with dotted boundaries). The path of the resonance frequency displays reverse shrinkage as its coefficient gets, in general, smaller for smaller values of the penalty. As expected from the theoretical results, if the penalty parameter is chosen correctly, it is possible to separate the signal and resonance frequencies for sufficiently low noise levels by just retaining large and neglecting small coefficients. It is also apparent that the coefficient of the resonance frequency is small for a correct choice of the penalty parameter but very seldom identically zero.*

The rate of convergence is actually optimal with an appropriate choice of the tuning parameter $\lambda$ and under the condition of bounded maximal and minimal sparse eigenvalues. This rate is, apart from logarithmic factor in $p_n$ and $n$, identical to what could be achieved if the true sparse model would be known. If $\ell_2$-consistency is achieved, the Lasso is selecting all "sufficiently large" coefficients, and possibly some other unwanted variables. "Sufficiently large" means here that the squared size of the coefficients is decaying slower than the rate $n^{-1} s_n \log p_n$, again ignoring logarithmic factors in the sample size. The number of variables can thus be narrowed down considerably with the Lasso in a meaningful way, keeping all important variables. The size of the reduced subset can be bounded with high probability by the number of truly important variables times a factor that depends on the decay of the sparse eigenvalues. This factor is often simply the squared logarithm of the sample size. Our conditions are similar in spirit to those in related aforementioned works, but expand the ground to cover possibly cases with more dependent predictors than UUP. These results support that the Lasso is a useful model identification method for high-dimensional data.

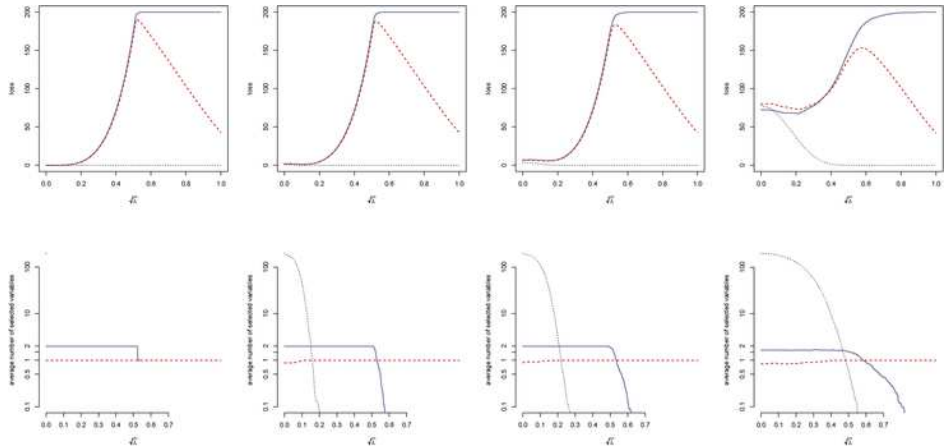

Fig. 3. *The top row shows the $\ell_2$-distance between $\beta$ and $\hat{\beta}^\lambda$ separately for the signal frequencies (solid blue line), resonance frequency (broken red line) and all other frequencies (dotted gray line). It is evident that the distance is quite small for all three categories simultaneously if the noise level is sufficiently low (the noise level is again increasing from left to right as $\sigma = 0, 0.1, 0.2, 1$). The bottom row shows, on the other hand, the average number of selected variables (with nonzero estimated regression coefficient) in each of the three categories as a function of the penalty parameter. It is impossible to choose the correct model, as the resonance frequency is always selected, no matter how low the noise level and no matter how the penalty parameter is chosen. This illustrates that sign consistency does not hold if the irrepresentable condition is violated, even though the estimate can be close to the true vector $\beta$ in the $\ell_2$-sense.*

**Acknowledgments.** We would like to thank Noureddine El Karoui and Debashis Paul for pointing out interesting connections to Random Matrix theory. Some results of this manuscript have been presented at the Oberwolfach workshop "Qualitative Assumptions and Regularization for High-Dimensional Data." Finally, we would like to thank the two referees and the AE for their helpful comments that have led to an improvement over our previous results.

Department of Statistics  
University of Oxford  
1 South Parks Road  
Oxford OX1 3TG  
United Kingdom  
E-mail: meinshausen@stats.ox.ac.uk

Department of Statistics  
UC Berkeley  
367 Evans Hall  
Berkeley, California 94720  
USA  
E-mail: yu@stat.berkeley.edu